\documentclass[12pt]{article}\def\today{May 30, 2000}
\usepackage{amssymb}
\usepackage{amsmath}
\usepackage[mathscr]{eucal}

\newtheorem{defA}{Definition}
\renewcommand{\thedefA}{\Alph{defA}}
\newenvironment{deA}{\begin{defA}\rm}{\end{defA}}

\newtheorem{thm}{Theorem}[section]

\newtheorem{lem}[thm]{Lemma}

\newtheorem{assumption}[thm]{Assumption}

\newtheorem{definition}[thm]{Definition}
\newenvironment{de}{\begin{definition}\rm}{\end{definition}}
\newtheorem{example}[thm]{Example}
\newenvironment{exmp}{\begin{example}\rm}{\end{example}}
\newtheorem{remark}[thm]{Remark}
\newenvironment{rem}{\begin{remark}\rm}{\end{remark}}
\newcommand{\A}{\mathcal{A}}
\newcommand{\B}{\mathcal{B}}
\newcommand{\C}{\mathcal{C}}
\newcommand{\F}{\mathbb{F}}

\newcommand{\Z}{\mathbb{Z}}

\newcommand{\Ham}{d_H}
\newcommand{\rank}{\mathrm{rank}\,}
\newcommand{\im}{\mathrm{im}\,}

\newcommand{\eqr}[1]{~\mbox{$(${\rm \ref{#1}}$)$}}
\newcommand{\Section}[1]{\section{#1}\setcounter{equation}{0}}

\newcommand{\openbox}{\leavevmode
  \hbox to.77778em{%
    \hfil\vrule
  \vbox to.675em{\hrule width.6em\vfil\hrule}%
  \vrule\hfil}} \newcommand{\proofname}{Proof}
\newenvironment{proof}[1][\proofname]{\par\normalfont
  \trivlist\item[\hskip\labelsep\itshape #1:]\ignorespaces
  }{\hspace*{1cm}\hspace*{\fill}\openbox \medskip\endtrivlist}

\newenvironment{proofs}[1][Sketch of Proof]{\par\normalfont
  \trivlist\item[\hskip\labelsep\itshape #1:]\ignorespaces
  }{\hspace*{1cm}\hspace*{\fill}\openbox \medskip\endtrivlist}

\title{Connections between Linear Systems and\\
  Convolutional Codes}%
\date{\today}%
\author{Joachim Rosenthal\footnote{Supported in part by NSF grant
    DMS-96-10389.  This research has been carried out while the
    author was a guest professor at EPFL in Switzerland. The
    author would like to thank EPFL for its support and
    hospitality.}\\
  {\small Department of Mathematics\vspace{-2mm}}\\
  {\small University of Notre Dame\vspace{-2mm}}\\
  {\small Notre Dame, Indiana 46556-5683, USA}\\
  {\small {\em e-mail:\/} Rosenthal.1@nd.edu} }

\begin{document} \maketitle
\begin{abstract}
  The article reviews different definitions for a convolutional
  code which can be found in the literature. The algebraic
  differences between the definitions are worked out in detail.
  It is shown that bi-infinite support systems are dual to
  finite-support systems under Pontryagin duality. In this
  duality the dual of a controllable system is observable and
  vice versa. Uncontrollability can occur only if there are
  bi-infinite support trajectories in the behavior, so finite and
  half-infinite-support systems must be controllable.
  Unobservability can occur only if there are finite support
  trajectories in the behavior, so bi-infinite and
  half-infinite-support systems must be observable.  It is shown
  that the different definitions for convolutional codes are
  equivalent if one restricts attention to controllable and
  observable codes.\bigskip

\end{abstract}

\Section{Introduction}

It is common knowledge that there is a close connection between
linear systems over finite fields and convolutional codes. In the
literature one finds however a multitude of definitions for
convolutional codes, which can make it confusing for somebody
who wants to enter this research field with a background in
systems theory or symbolic dynamics.  It is the purpose of this
article to provide a survey of the different points of view about
convolutional codes.

The article is structured as follow: In Section~\ref{Lin-Alg} we
will review the way convolutional codes have often been defined
in the coding literature~\cite{jo93,jo99,li83,mc98,pi88}.

Section~\ref{Sym-Dyn} reviews a definition for convolutional
codes that can be found in the literature on symbolic dynamics.
{}From the symbolic dynamics point of view~\cite{ki00u,li95,ma95},
a convolutional code is a linear irreducible shift space.

In Section~\ref{Wil-Beh} we will review the class of
time-invariant, complete linear behaviors in the sense of
Willems~\cite{wi86a1,wi89,wi91}. We will show how these behaviors
relate to the definitions given in Section~\ref{Lin-Alg}
and~\ref{Sym-Dyn}.

In Section~\ref{Module} we will give a definition for
convolutional codes in which it is required that the code words
have finite support. Such a definition was considered by
Fornasini and Valcher~\cite{va94,fo98a} and by the author in
collaboration with Schumacher, Weiner and
York~\cite{ro96a1,ro99a,we98t}. The study of behaviors with
finite support has been done earlier in the context of automata
theory and we refer to Eilenberg's book~\cite{ei74}.  We show in
Section~\ref{Module} how this module-theoretic definition relates
to complete, linear and time-invariant behaviors by Pontryagin
duality.

In Section~\ref{First-Order} we will study different first-order
representations connected with the different viewpoints.  Finally,
in Section~\ref{Case1} we compare the different definitions.  We
also show how cyclic redundancy check codes can naturally be
viewed in the context of finite-support convolutional codes.

Throughout the paper we will emphasize the algebraic properties
of the different definitions. We will also restrict ourselves to
the concrete setting of convolutional codes defined over finite
fields. It is however known that many of the concepts in this
paper generalize to group codes~\cite{fa99,fo00u,fo99} and
multidimensional convolutional
codes~\cite{fo95p2,fo98a,gl99u,va94,we98t}.  All of the definitions
which we are going to give are quite similar, but there are some
notable differences.

Since the paper draws from results from quite different research
areas, one is faced with the problem that there is no uniform
notation. In this paper we will adopt the convention used in
systems theory in which vectors are regarded as column vectors.
For the convenience of the reader, we conclude this section with a
summary of some of the notation used in this paper:

\begin{tabular}{ll}
$\F$ & A fixed finite field;\\
$\F[z]$ & The polynomial ring over $\F$;\\
$\F[z,z^{-1}]$ & The Laurent polynomial ring over $\F$;\\
$\F(z)$ & The field of rationals;\\
$\F[[z]]$ & The ring of formal power  series of the form
$\sum_{i=0}^\infty a_iz^i$;\\
$\F((z))$ & The field of formal Laurent series having the form
$\sum_{i=d}^\infty a_iz^i$; \\
$\F[[z,z^{-1}]]$ & The ring of formal power series of the form
$\sum_{i=-\infty}^\infty a_iz^i$;\\
$\Z$ & The integers;\\
$\Z_+$ & The nonnegative integers;\\
$\Z_-$ & The nonpositive integers.
\end{tabular}

Consider the ring of formal power series $\F[[z,z^{-1}]]$. We
will identify the set $\F[[z,z^{-1}]]$ with the (two-sided)
sequence space $\F^\Z$.  We have natural embeddings:
$$
\F\longrightarrow \F[z]\longrightarrow \F[z,z^{-1}]
\longrightarrow \F(z)\longrightarrow \F((z))\longrightarrow
\F[[z,z^{-1}]].
$$
With these embeddings we can view e.g. the set of rationals
$\F(z)$ as a subset of the sequence space $\F^\Z$, and we will
make use of such identifications throughout the paper.

The set of $n$-vectors with polynomial entries will be denoted by
$\F^n[z]$. Similarly we define the sets $\F^n(z),\F^n((z))$ etc.
All these sets are subsets of the two sided sequence space
$\left( \F^n\right)^\Z=\F^n[[z,z^{-1}]]$. The definitions of
convolutional codes which we will provide in the next sections
will all be $\F$-linear subspaces of $\left( \F^n\right)^\Z$.

The idea of writing a survey on the different points of view
about convolutional codes was suggested to the author by Paul
Fuhrmann during a stimulating workshop on ``Codes, Systems and
Graphical Models'' at the Institute for Mathematics and its
Applications (IMA) in August 1999.  A first draft of this paper
was circulated in October 1999 to about a dozen people interested
in these research issues. This generated an interesting `Internet
discussion' on these issues, in which the different opinions were
exchanged by e-mail. Some of these ideas have been incorporated
into the final version of the paper and the author would like to
thank Dave Forney, Paul Fuhrmann, Heide Gluesing-Luerssen, Jan
Willems and Sandro Zampieri for having provided valuable
thoughts. The author wishes also to thank the IMA and its superb
staff, who made the above mentioned workshop possible.

\Section{The linear algebra point of view} \label{Lin-Alg}

The theory of convolutional codes grew out and extended the
theory of linear block codes into a new direction. Because of
this reason we start the section with linear block codes and we
introduce convolutional codes in a quite intuitive way.

An $[n,k]$ linear block code is by definition a linear subspace
$\C\subset\F^n$ having dimension $\dim\C=k$. Let $G$ be a
$n\times k$ matrix with entries in $\F$. The linear map
$$
\varphi:\ \F^k\longrightarrow\F^n,\ m\longmapsto c=Gm
$$
is called an {\em encoding map\/} for the code $\C$ if
$\im(\varphi)=\C$.  If this is the case then we say $G$ is a {\em
  generator matrix\/} or an {\em encoder\/} for the block code
$\C$.

Assume that a sequence of message blocks
$m_0,\ldots,m_t\subset\F^k$ should be encoded into a
corresponding sequence of code words $c_i=Gm_i\in\F^n,\ 
i=0,\ldots,t$. By introducing the polynomial vectors
$m(z)=\sum_{i=0}^tm_iz^i\in{\F}^k[z]$ and
$c(z)=\sum_{i=0}^tc_iz^i\in{\F}^n[z]$ it is possible to describe
the encoding procedure through the module
homomorphism:\footnote{Throughout the paper we use the symbol
  $\varphi$ to denote an encoding map. The context will make it
  clear what the domain and the range of this map is in each
  situation.}
\begin{equation}    \label{enco2}
\varphi:\ \F^k[z]\longrightarrow\F^n[z],\ m(z)\longmapsto
c(z)=Gm(z). 
\end{equation}

The original idea of a convolutional code goes back to the paper
of Elias~\cite{el55}, where it was suggested to use a polynomial
matrix $G(z)$ in the encoding procedure\eqr{enco2}.

Polynomial encoders $G(z)$ are physically easily implemented
through a feedforward linear sequential circuit. Massey and
Sain~\cite{ma67,sa69} showed that there is a close connection
between linear systems and convolutional codes. Massey and Sain
viewed the polynomial encoder $G(z)$ as a transfer function.
More generally it is possible to realize a transfer function
$G(z)$ with rational entries by (see e.g.~\cite{jo93,jo99}) a
linear sequential circuit whose elements include feedback
components. If one allows rational entries in the encoding matrix
then it seems natural to extend the possible message sequences to
the set of rational vectors $m(z)\in\F^k(z)$ and to process this
sequence by a `rational encoder' resulting again in a rational
code vector $c(z)\in\F^n(z)$. With this we have a first
definition of a convolutional code as it can be found e.g. in the
Handbook of Coding Theory~\cite[Definition~2.4]{mc98}:
\begin{deA}       \label{Def-A}
  A $\F(z)$-linear subspace $\C$ of $\F^n(z)$ is called a
  convolutional code.
\end{deA}

If $G(z)$ is a $n\times k$ matrix with entries in $\F(z)$ whose
columns form a basis for $\C$, then we call $G(z)$ a generator
matrix or an encoder for the convolutional code $\C$. $G(z)$
describes the encoding map:
$$
\varphi:\ \F^k(z)\longrightarrow\F^n(z),\ m(z)\longmapsto
c(z)=G(z)m(z).
$$
The field of rationals $\F(z)$ viewed as a subset of the
sequence space $\F^\Z=\F[[z,z^{-1}]]$ consists precisely of those
sequences whose support is finite on the negative sequence space
$\F^{\Z_{-}}$ and whose elements form an ultimately periodic
sequence on the positive sequence space $\F^{\Z_{+}}$. It
therefore seems that one equally well could restrict the possible
message words $m(z)\in\F^k(z)$ to sequences whose coordinates
consists of Laurent polynomials only, in other words to sequences
of the form $m(z)\in\F^k[z,z^{-1}]$.

Alternatively one could allow message words $m(z)$ whose
coordinates are not ultimately periodic and possibly not of
finite support on the negative sequence space $\F^{\Z_{-}}$. This
would suggest that one should take as possible message words the
whole sequence space $\left(\F^k\right)^\Z=\F^k[[z,z^{-1}]]$. The
problem with this approach is that the multiplication of an
element in $\F[[z,z^{-1}]]$ with an element in $\F(z)$ is in
general not well defined. If one restricts however the message
sequences to the field of formal Laurent series then the
multiplication is well defined. This leads to the following
definition which goes back to the work of Forney~\cite{fo73}.
The definition has been adopted in the book by Piret~\cite{pi88}
and the book by Johannesson and Zigangirov~\cite{jo99}, and it
appears as Definition~2.3 in the Handbook of Coding
Theory~\cite{mc98}:

\renewcommand{\thedefA}{\Alph{defA}$^\prime$}  \addtocounter{defA}{-1}
\begin{deA}                   \label{Def-A-pr}
  A $\F((z))$-linear subspace $\C$ of $\F^n((z))$ which has a
  basis of rational vectors in $\F^n(z)$ is called a
  convolutional code.
\end{deA}
\renewcommand{\thedefA}{\Alph{defA}}

The requirement that $\C$ has a basis with rational entries
guarantees that $\C$ has also a basis with only polynomial
entries. $\C$ can therefore be represented by a $n\times k$
generator matrix $G(z)$ whose entries consist only of rationals
or even even polynomials. The encoding map with respect to $G(z)$
is given through:
\begin{equation}    \label{map-A}
\varphi:\ \F^k((z))\longrightarrow\F^n((z)),\ m(z)\longmapsto
c(z)=G(z)m(z).
\end{equation}

If $G(z)$ is a polynomial matrix, then finitely many components of
$m(z)$ influence only finitely many components of $c(z)$, and the
encoding procedure may be physically implemented by a simple
feedforward linear shift register.

If $G(z)$ contains rational entries, then it is in general the
case that a finite (polynomial) message vector is encoded into an
infinite (rational) code vector of the form
$c(z)=\sum_{i=s}^\infty c_iz^i$. This might cause some
difficulties in the decoder. For the encoding process, $G(z)$ can
be physically realized by linear shift registers, in general with
feedback (see e.g.~\cite{jo93,jo99}).

{}From a systems theory point of view, it is
classical~\cite{ka69} to view the encoding map\eqr{map-A} as an
input-output linear system. This was the point of view taken by
Massey and Sain~\cite{ma67,sa69} and thereafter in most of the
coding literature. However unlike in systems theory, the
important object in coding theory is the code $\C=\im(\varphi)$.
As a result one calls encoders $\varphi$ which generate the same
image $\im(\varphi)$ equivalent; we will say more about this in a
moment. In Sections~\ref{Sym-Dyn} and~\ref{Wil-Beh} we will
view\eqr{map-A} as an image representation of a time-invariant
behavior in the sense of Willems~\cite{wi86a1,wi89}, which we
believe captures the coding situation in a more natural way.

Assume that $G(z)$ and $\tilde{G}(z)$ are two $n\times k$
rational encoding matrices defining the same code $\C$ with
respect to either Definition~\ref{Def-A} or~\ref{Def-A-pr}. In
this case we say that $G(z)$ and $\tilde{G}(z)$ are equivalent
encoders. The following lemma is a simple result of linear
algebra:
\begin{lem}                      \label{lem-eq-A}
  Two $n\times k$ rational encoders $G(z)$ and $\tilde{G}(z)$ are
  equivalent with respect to either Definition~\ref{Def-A}
  or~\ref{Def-A-pr} if and only if there is a $k\times k$
  invertible rational matrix $R(z)$ such that
  $\tilde{G}(z)=G(z)R(z)$.
\end{lem}

It follows from this lemma that Definition~\ref{Def-A} and
Definition~\ref{Def-A-pr} are completely equivalent with respect
to equivalence of encoders.

{}From an algebraic point of view we can identify a convolutional
code in the sense of Definition~\ref{Def-A} or
Definition~\ref{Def-A-pr} through an equivalence class of
rational matrices. The following theorem singles out a set of
very desirable encoders inside each equivalence class.
\begin{thm}                      \label{basicM}
  Let $G(z)$ be a $n\times k$ rational encoding matrix of
  rank~$k$ defining a code $\C$. Then there is a $k\times k$
  invertible rational matrix $R(z)$ such that
  $\tilde{G}(z)=G(z)R(z)$ has the properties:
  \begin{description}
  \item[(i)] $\tilde{G}(z)$ is a polynomial matrix.
  \item[(ii)] $\tilde{G}(z)$ is right prime.
  \item[(iii)] $\tilde{G}(z)$ is column reduced with column
    degrees $\{e_1, \ldots, e_k\}$.
  \end{description}
  Furthermore, every polynomial encoding matrix of $\C$ which is
  right prime and column-reduced has (unordered) column degrees
  $\{e_1, \ldots, e_k\}$. Thus these indices are invariants of
  the convolutional code.
\end{thm}

The essence of Theorem~\ref{basicM} was proved by
Forney~\cite[Theorem 3]{fo70}. In~\cite{fo75} Forney related the
indices appearing in (iii) to the controllability and
observability indices of a controllable and observable system.
Paper~\cite{fo75} had an immense impact in the linear systems
theory literature. We will follow here the suggestion of
McEliece~\cite{mc98} and call these indices the {\em Forney
  indices\/} of the convolutional code, despite the fact that
Theorem~\ref{basicM} can be traced back to the last century, when
Kronecker, Hermite and in particular Dedekind and Weber studied
matrices over the rationals and more general function fields.  In
Sections~\ref{Wil-Beh} and~\ref{Module} we will make a
distinction between the Forney indices as defined above and the
Kronecker indices of a submodule of $\F^n[z]$.

In the coding literature~\cite{jo99,pi88}, an encoder satisfying
conditions (i), (ii) and (iii) of Theorem~\ref{basicM} is called
a {\em minimal basic encoder\/}.

So far we have used encoding matrices to describe a convolutional
code. As is customary in linear algebra, one often describes a
linear subspace as the kernel of a matrix. This leads to the
notion of a {\em parity-check matrix}. The following theorem is
well known (see e.g.~\cite{pi88}).
\begin{thm}                      \label{basicPar}
  Let $\C\subset \F^n((z))$ be a rank-$k$ convolutional code in
  the sense of Definition~\ref{Def-A-pr}. Then there exists an
  $r\times n$ matrix $H(z)$ such that the code is equivalently
  described as the kernel of $H(z)$:
  $$
  \C=\{ \ c(z)\in\F^n((z))\ \mid\ \ H(z)c(z)=0\ \}.
  $$
  Moreover, it is possible to choose $H(z)$ in such a way that:
  \begin{description}
  \item[(i)] $H(z)$ is a polynomial matrix.
  \item[(ii)] $H(z)$ is left prime.
  \item[(iii)] $H(z)$ is row-reduced having row
    degrees $\{f_1, \ldots, f_r\}$.
  \end{description}
  Furthermore, every polynomial parity check matrix of $\C$ which
  is left prime and row reduced will have (unordered) row degrees
  $\{f_1, \ldots, f_r\}$. Thus these indices are invariants of
  the convolutional code.
\end{thm}

Properties (i)--(iii) essentially follow from the fact that the
transpose $H^t(z)$ is a generator matrix for the dual
(orthogonal) code $\C^\perp$.

The set of indices $\{e_1, \ldots, e_k\}$ and $\{f_1, \ldots,
f_r\}$ differ in general, their sum is however always the same,
and is called the {\em degree\/} of the convolutional code.  One
says that a rank-$k$ code $\C\subset \F^n((z))$ has {\em
  transmission rate\/} $k/n$, {\em controller memory\/}
$m:=\max\{e_1, \ldots, e_k\}$ and {\em observer memory\/}
$n:=\max \{f_1, \ldots, f_r\}$.
  
Another important code parameter is the {\em free distance}.
The free distance of a code measures the smallest distance
between any two different code words, and is formally defined
as:
\begin{equation}                     \label{free-dist}        
 d_\mathrm{free}(\C):=\min_{{u,v \in\C\atop u\neq v}} 
 \sum_{t\in \Z}
 \Ham(u_t,v_t),
\end{equation}
where $\Ham(\ , \ )$ denotes the usual Hamming distance on $\F^n$.

\Section{The symbolic dynamics point of view}
\label{Sym-Dyn}

In this section we present a definition of convolutional codes as
it can be found in the symbolic dynamics
literature~\cite{ki00u,li95,ma95}. Convolutional codes in this
framework are exactly the linear, compact, irreducible and
shift-invariant subsets of $\F^n[[z,z^{-1}]]$. In order to make
this precise, we will have to develop some basic notions from
symbolic dynamics.

In the sequel we will work with the finite alphabet $\A:=\F^n$. A
{\em block\/} over the alphabet $\A$ is a finite sequence
$\beta=x_1 x_2\ldots x_k$ consisting of $k$ elements $x_i\in\A$.
If $w=w(z)=\sum_iw_iz^i\in\F^n[[z,z^{-1}]]$ is a sequence, one
says that the block $\beta$ occurs in $w$ if there is some
integer $j$ such that $\beta=w_jw_{j+1}\ldots w_{k+j-1}$. If
$X\subset \F^n[[z,z^{-1}]]$ is any subset, we denote by
$\mathscr{B}(X)$ the set of blocks which occur in some element of
$X$.

The fundamental objects in symbolic dynamics are the {\em
  shift spaces.\/} For this let $\mathscr{F}$ be a set of blocks,
possibly infinite.
\begin{de}
  The subset $X\subset \F^n[[z,z^{-1}]]$ consisting of all
  sequences $w(z)$ which do not contain any of the (forbidden)
  blocks of $\mathscr{F}$ is called a {\em shift space}.
\end{de}

The left-shift operator is the $\F$-linear map
\begin{equation}                \label{leftshift}
\sigma:\ \F[[z,z^{-1}]]\longrightarrow \F[[z,z^{-1}]],\ \ 
w(z)\longmapsto z^{-1}w(z).
\end{equation}
Let $I_n$ be the $n\times n$ identity matrix. The shift map
$\sigma$ extends to the shift map
$$
\sigma I_n:\ \F^n[[z,z^{-1}]]\longrightarrow \F^n[[z,z^{-1}]].
$$
One says that $X\subset \F^n[[z,z^{-1}]]$ is a {\em shift-invariant
  set\/} if $(\sigma I_n)(X)\subset X$. Clearly shift spaces are
shift-invariant subsets of $\F^n[[z,z^{-1}]]$.

It is possible to characterize shift spaces in a topological
manner. For this we will introduce a metric on
$\F^n[[z,z^{-1}]]$:
\begin{de}                   \label{defmetric}
  If $v(z)=\sum_iv_iz^i$ and $w(z)=\sum_iw_iz^i$ are both
  elements of $\F^n[[z,z^{-1}]]$ we define their distance
  through:
  \begin{equation}                   \label{def2metric}
    d(v(z),w(z)):=\sum_{i\in\Z}2^{-|i|}\Ham(v_i,w_i).
  \end{equation}
\end{de}
In this metric two elements $v(z),w(z)$ are `close' if they
coincide over a `large block around zero'. One readily verifies
that $d(\ , \ )$ indeed satisfies all the properties of a metric
and therefore induces a topology on $\F^n[[z,z^{-1}]]$. Using
this topology we can characterize shift spaces:
\begin{thm}
  A subset of $\F^n[[z,z^{-1}]]$ is a shift space if and only if
  it is shift-invariant and compact.
\end{thm}
\begin{proof}
  The metric introduced in Definition~\ref{defmetric} is
  equivalent to the metric described
  in~\cite[Example~6.1.10]{li95}. The induced topologies are
  therefore the same. The result follows therefore
  from~\cite[Theorem~6.1.21]{li95}.
\end{proof}

The topological space $\F^n[[z,z^{-1}]]$ is a typical example of
a linearly compact vector space, a notion introduced by
S.~Lefschetz. There is a large theory on linearly compact vector
spaces, and several of the results which we are going to derive
are valid in this broader context. We refer the interested reader
to~\cite[\S 10]{ko69} for more details.

A further important concept is irreducibility which will turn out
to be equivalent to the concept of controllability in our
concrete setting.
\begin{de}                            \label{irreducible}
  A shift space $X\subset \F^n[[z,z^{-1}]]$ is called {\em
    irreducible\/} if for every ordered pair of blocks
  $\beta,\gamma$ of $\mathscr{B}(X)$ there is a block $\mu$ such
  that the concatenated block $\beta\mu\gamma$ is in
  $\mathscr{B}(X)$.
\end{de}

We are now prepared to give the symbolic dynamics definition for
a convolutional code and to work out the basic properties for
these codes.

\begin{deA}                   \label{Def-B}
  A linear, compact, irreducible and shift-invariant subset of
  $\F^n[[z,z^{-1}]]$ is called a convolutional code.
\end{deA}

This is an abstract definition and it is not immediately clear
how one should encode messages with such convolutional codes.
The following will make this clear.

Let $G(z)$ be a $n\times k$ matrix with entries in the ring of
Laurent polynomials $\F[z,z^{-1}]$. Consider the encoding map:

\begin{equation}              \label{map-B}  
\varphi:\ \F^k[[z,z^{-1}]]\longrightarrow\F^n[[z,z^{-1}]], \ 
m(z)\longmapsto c(z)=G(\sigma)m(z).
\end{equation}
In terms of polynomials the map $\varphi$ is simply described
through $m(z)\longmapsto c(z)=G(z^{-1})m(z)$.

Recall that a continuous map is called closed if the image of a
closed set is closed. Using the fact that $\F^n[[z,z^{-1}]]$ is
compact, one (easily) proves the following result:
\begin{lem}                            \label{cont+closed}
  The encoding map\eqr{map-B} is $\F$-linear, continuous and
  closed.
\end{lem}
Clearly $\im(\varphi)$ is also shift-invariant, and one
shows~\cite{li95} that the image of an irreducible set under
$\varphi$ is irreducible again.
 
In summary we have shown that $\im(\varphi)$ describes a
convolutional code in the sense of Definition~\ref{Def-B}.
Actually the converse is true as well:
\begin{thm}                       \label{image}
  $\C\subset\F^n[[z,z^{-1}]]$ is a convolutional code in the
  sense of Definition~\ref{Def-B} if and only if there exists a
  Laurent polynomial matrix $G(z)$ such that $\C=\im(\varphi)$,
  where $\varphi$ is the map in\eqr{map-B}.
\end{thm}

A proof of this theorem will be given in the next section after
Theorem~\ref{wi2}.

The question now arises how Definition~\ref{Def-B} relates to
Definition~\ref{Def-A} and Definition~\ref{Def-A-pr}.  The
following theorem will provide a partial answer to this question.
\begin{thm}                      \label{closure}
  Assume that $\C\subset\F^n[[z,z^{-1}]]$ is a nonzero
  convolutional code in the sense of Definition~\ref{Def-A} or
  Definition~\ref{Def-A-pr}. Then $\C$ is not closed, but the
  closure of $\bar{\C}$ of $\C$ is a convolutional code in the
  sense of Definition~\ref{Def-B}.
\end{thm}
\begin{proof}
  Let $G(z)$ be a minimal basic encoder of $\C$ and let
  $w(z)\in\F^n[z]$ be the first column of $G(z)$. Note that $w(z)
  \in\C$ and that there is at least one entry of $w(z)$ which
  does not contain the factor $(z-1)$. Let
  $\phi_N(z):=\sum_{i=-N}^Nz^i\in\F[z,z^{-1}]$ and consider the
  sequence of code words $w^N(z):=\phi_N(z)w(z)$.  For each $N>0$
  one has that $w^N(z)\in\C$.  However
  $\lim_{N\mapsto\infty}w^N(z)$ is in $\F^n[[z,z^{-1}]]\setminus
  \F^n((z))\subset \F^n[[z,z^{-1}]]\setminus\C$. This shows that
  $\C$ is not a closed set inside $\F^n[[z,z^{-1}]]$. The closure
  $\bar{\C}$ is obtained by extending the input space $F^k((z))$
  to all of $F^k[[z,z^{-1}]]$. The image of $F^k[[z,z^{-1}]]$
  under the encoding map\eqr{map-B} is closed by
  Lemma~\ref{cont+closed}, hence the closure is a code in the
  sense of Definition~\ref{Def-B}.
\end{proof}

Actually one can show that there is a bijective correspondence
between the convolutional codes in the sense of
Definition~\ref{Def-A} (respectively Definition~\ref{Def-A-pr})
and the convolutional codes in the sense of
Definition~\ref{Def-B}, as we will show in~Theorem~\ref{Thm-7-1}
and~Theorem~\ref{Thm-7-2}. It is also worthwhile to remark that
already in 1983 Staiger published a paper~\cite{st83} where he
studied the closure of convolutional codes generated by a
polynomial generator matrix.

In analogy to Lemma~\ref{lem-eq-A}, one has:
\begin{lem}                      \label{lem-eq-B}
  Two $n\times k$ encoding matrices $G(z)$ and $\tilde{G}(z)$
  defined over the Laurent polynomial ring $\F[z,z^{-1}]$ are
  equivalent with respect to Definition~\ref{Def-B} if and only
  if there is a $k\times k$ invertible rational matrix $R(z)$
  such that $\tilde{G}(z)=G(z)R(z)$.
\end{lem}

We leave the proof again as an exercise for the reader. We remark
that rational transformations of the form $R(z)$ are needed to
describe the equivalence, even though it is in general not
possible to use a rational encoder $G(z)$ in the encoding
procedure\eqr{map-B}. This is simply due to the fact that in
general the multiplication of an element of $\F(z)$ with an
element of $\F[[z,z^{-1}]]$ is not defined. The following example
should make this clear. (Compare also with
Remark~\ref{rational}.)
\begin{exmp}                   \label{not-def}
  Consider $f(z)=\frac{1}{1-z}=\sum_{i=0}^\infty z^i\in\F(z)$ and
  $g(z)=\sum_{i=-\infty}^\infty z^i\in\F[[z,z^{-1}]]$. Trying to
  multiply the two power series $f(z),g(z)$ would result in a
  power series in which each coefficient would be infinite.
\end{exmp}

In the same way as at the end of Section~\ref{Lin-Alg} we define
the transmission rate, the degree, the memory and the free
distance of a convolutional code $\C$ in the sense of
Definition~\ref{Def-B}.

\Section{Linear time-invariant behaviors}
\label{Wil-Beh}

In this section we will take the point of view that a
convolutional code is a linear time-invariant behavior in the
sense of Willems~\cite{wi86a1,wi89,wi91}. Of course behavioral
system theory is quite general, allowing all kinds of time axes
and signal spaces. In order to relate the behavioral concepts to
the previous points of view, we will restrict our study to linear
behaviors in $\left(\F^n\right)^\Z=\F^n[[z,z^{-1}]]$ and
$\left(\F^n\right)^{\Z_+}=\F^n[[z]]$.

Let $\sigma$ be the shift operator defined in\eqr{leftshift}. One
says that a subset $\B\subset \F^n[[z,z^{-1}]]$ is {\em
  time-invariant\/} if $(\sigma I_n)(\B)\subset \B$. The concept
therefore coincides with the symbolic dynamics concept of
shift-invariance.

In addition to linearity and time-invariance, there is a third
important concept usually required of a time-invariant behavior:

\begin{de}                          \label{defcomplete}
  A behavior $\B\subset \F^n[[z,z^{-1}]]$ is said to be {\em
    complete\/} if $w \in \F^n[[z,z^{-1}]]$ belongs to $\B$
  whenever $w|_J$ belongs to $\B|_J$ for every finite subinterval
  $J\subset \Z$.
\end{de}

The definition simply says that $\B$ is complete if membership
can be decided on the basis of finite windows.  Completeness is
an important well behavedness property for linear time-invariant
behaviors, as Willems~\cite[p. 567]{wi86a1} emphasized with the
remark:
\begin{quote}
  As such, it can be said that the study of non-complete systems
  does not fall within the competence of system theorists and
  could be left to cosmologists or theologians.
\end{quote}

In Definition~\ref{defmetric} we introduced a metric on the
vector space $\F^n[[z,z^{-1}]]$. We remark that with respect to
this metric a subset $\B\subset \F^n[[z,z^{-1}]]$ is complete if
and only if every Cauchy sequence converges inside $\B$. In other
words, the completeness notion of Definition~\ref{defcomplete}
coincides with the usual topological notion of completeness.

The following result is known for linearly compact vector spaces,
a proof can be found in~\cite{wi86a1}:
\begin{lem}                         \label{closedcomplete}
  A linear subset $\B\subset\F^n[[z,z^{-1}]]$ is complete if and
  only if it is closed and hence compact.
\end{lem}
With these preliminaries we can define a convolutional code as
follows:
\begin{deA}                         \label{Def-C}
  A linear, time-invariant and complete subset
  $\B\subset\F^n[[z,z^{-1}]]$ is called a convolutional code.
\end{deA}
It is immediate from Lemma~\ref{closedcomplete} that the
convolutional codes defined in Definition~\ref{Def-B} are
complete and that Definition~\ref{Def-C} is more general than
Definition~\ref{Def-B}, since no irreducibility is required. It
also follows from Theorem~\ref{closure} and
Lemma~\ref{closedcomplete} that the convolutional codes defined
in Definition~\ref{Def-A} and Definition~\ref{Def-A-pr} are in
general not complete.

Before we elaborate on these differences we would like also to
treat the situation when the time axis is $\Z_+$ since
traditionally a large part of linear systems theory has been
concerned with systems defined on the positive time axis.

We first define the left-shift operator acting on
$\left(\F^n\right)^{\Z_+}=\F^n[[z]]$ through:
\begin{equation}                \label{leftshift+}
\sigma:\ \F[[z]]\longrightarrow \F[[z,z^{-1}]],\ \ 
w(z)\longmapsto z^{-1}(w(z)-w(0)).
\end{equation}
We have used the same symbol as in\eqr{leftshift} since the
context will always make it clear if we work over $\Z$ or $\Z_+$.
In analogy to\eqr{leftshift} $\sigma$ extends to the shift map $
\sigma I_n:\ \F^n[[z]]\longrightarrow \F^n[[z]]$, and one says a
subset $X\subset \F^n[[z]]$ is time-invariant if $(\sigma
I_n)(X)\subset X$. Notice however that the map of\eqr{leftshift+}, 
unlike that of\eqr{leftshift}, is not invertible.

With this we have: \renewcommand{\thedefA}{\Alph{defA}$^\prime$}
\addtocounter{defA}{-1}
\begin{deA}                          \label{Def-C-pr}
  A linear, time-invariant and complete subset
  $\B\subset\F^n[[z]]$ is called a convolutional code.
\end{deA}
\renewcommand{\thedefA}{\Alph{defA}}

The following fundamental theorem was proved by
Willems~\cite[Theorem 5]{wi86a1}.
\begin{thm}                          \label{thm-Wil}
  A subset $\B\subset\F^n[[z,z^{-1}]]$ (respectively a subset
  $\B\subset\F^n[[z]]$) is linear, time-invariant and complete if
  and only if there is a $r\times n$ matrix $P(z)$ having entries
  in $\F[z]$ such that
 \begin{equation}                        \label{kernelbehav}
   \B=\left\{ \ w(z)\ \mid\  P(\sigma)w(z)=0 \ \right\}.
 \end{equation}
\end{thm}
By Lemma~\ref{cont+closed} the linear map $\psi:\ 
\F^n[[z,z^{-1}]]\longrightarrow \F^n[[z,z^{-1}]],\ 
w(z)\longmapsto P(\sigma)w(z)$ is continuous and its kernel is
therefore a complete set. It is therefore immediate that the
behavior defined in\eqr{kernelbehav} is linear, time-invariant and
complete.  The harder part of Theorem~\ref{thm-Wil} is the
converse statement.

Equation\eqr{kernelbehav} is often referred to as a kernel (or
AR) representation of a behavioral system. We will denote a
behavior having the form\eqr{kernelbehav} by $\ker P(\sigma)$.
By contrast, the encoding map $\varphi$ defined in\eqr{map-B}
describes an image (or MA) representation of the behavior
$\im(\varphi)=\im G(\sigma)$.

The most general representation is an ARMA representation. For
this let $P(z)$ and $G(z)$ be matrices of size $r\times n$ and
$r\times k$ respectively, having entries in the Laurent
polynomial ring $\F[z,z^{-1}]$. Then
\begin{equation}                        \label{armabehav}
   \B=\left\{\ w(z)\in \F^n[[z,z^{-1}]]  \ \mid \ \exists
     m(z)\in\F^k[[z,z^{-1}]]:\ \
     P(\sigma)w(z)=G(\sigma)m(z) \ \right\}
\end{equation}
is called an ARMA model.  One immediately verifies that the set
$\B$ is linear and time-invariant. It is a direct consequence
of Lemma~\ref{cont+closed} that $\B$ is also closed and hence
complete. Theorem~\ref{thm-Wil} therefore states that it is
possible to eliminate the so called `latent variable' $m(z)$ and
describe the behavior $\B$ by a simpler kernel representation of
the form\eqr{kernelbehav}. It follows in particular that the code
$\im(\varphi)=\im G(\sigma)$ defined in\eqr{map-B} has an
equivalent kernel representation of the form\eqr{kernelbehav} but
that in general the converse is not true.

\begin{rem}              \label{rational}
  As we explained in Section~\ref{Lin-Alg} it is quite common to
  use rational encoders for convolutional codes. In the ARMA
  model\eqr{armabehav} we required that the entries of $P(z)$ and
  $G(z)$ be from the Laurent polynomial ring. If $P(z)$ and
  $G(z)$ were rational matrices, then the behavior $\B\subset
  \F^n[[z,z^{-1}]]$ appearing in\eqr{armabehav} might not be well
  defined, as we showed in Example~\ref{not-def}. On the other
  hand if one restricts the behavior to the positive time axis
  $\Z_+$, i.e.  if one assumes that $\B\subset \F^n[[z]]$, then
  the set\eqr{armabehav} is defined even if $P(z)$ and $G(z)$ are
  rational encoders. This is certainly one reason why much
  classical system theory focused on shift spaces $\B\subset
  \F^n[[z]]$ or $\B\subset \F^n((z))$.
\end{rem}

In the sequel we will concentrate on representations of the
form\eqr{kernelbehav}.  Again the question arises, when are two
kernel representations equivalent?
\begin{lem}                      \label{lem-eq-C}
  Two $r\times n$ matrices $P(z)$ and $\tilde{P}(z)$ defined over
  the Laurent polynomial ring $\F[z,z^{-1}]$ describe the same
  behavior $\ker P(\sigma)=\ker \tilde{P}(\sigma)\subset
  \F^n[[z,z^{-1}]]$ if and only if there is a $r\times r$ matrix
  $U(z)$, unimodular over $\F[z,z^{-1}]$, such that
  $\tilde{P}(z)=U(z)P(z)$.
\end{lem}
\begin{proof}
  \cite[Proposition III.3]{wi91}.
\end{proof}
Similarly, if $P(z)$ and $\tilde{P}(z)$ are defined over $\F[z]$,
then these matrices define the same behavior $\ker P(\sigma)=\ker
P(\sigma)\subset \F^n[[z]]$ if and only if there is a matrix
$U(z)$, unimodular over $\F[z]$, such that
$\tilde{P}(z)=U(z)P(z)$.

The major difference between Definition~\ref{Def-B} and
Definition~\ref{Def-C} seems to be that Definition~\ref{Def-C}
does not require irreducibility. This last concept corresponds to
the term controllability (see~\cite{fo95}) in systems theory. We
first start with some notation taken from~\cite{ro96a1}:
 
For a sequence
$w=\sum_{-\infty}^{\infty}w_iz^i\in\F^n[[z,z^{-1}]]$, we use the
symbol $w^+$ to denote the `right half' $\sum_0^\infty w_iz^i$
and the symbol $w^-$ to denote the `left half' $\sum_{-\infty}^0
w_iz^i$.

\begin{de}                        \label{control}
  A behavior $\B$ defined on $\Z$ is said to be {\em
    controllable\/} if there is some integer $\ell$ such that for
  every $w$ and $w'$ in $\B$ and every integer $j$ there exists a
  $w'' \in \B$ such that $(z^{j}w'')^- = (z^{j}w)^-$ and
  $(z^{j+\ell}w'')^+ = (z^{j+\ell}w')^+$.
\end{de}

\begin{rem}
  Loeliger and Mittelholzer~\cite{lo96a} speak of {\em strongly
    controllable} if a behavior satisfies the conditions of
  Definition~\ref{control}. `Weakly controllable' in contrast
  requires an integer~$\ell$ which may depend on the trajectories
  $w$ and $w'$. The notions are equivalent in our concrete
  setting.
\end{rem}

We leave it as an exercise for the reader to show that
irreducibility as introduced in Definition~\ref{irreducible} is
equivalent to controllability for linear, time-invariant and
complete behaviors $\B\subset\F^n[[z,z^{-1}]]$. The next theorem
gives equivalent conditions for a behavior to be controllable.

\begin{thm}                               \label{wi2}
  (cf.\,\cite[Prop.\,4.3]{wi89}) Let $P(z)$ be a $r\times n$
  matrix of rank $r$ defined over $\F[z,z^{-1}]$. The following
  conditions are equivalent:
\begin{description}
\item[(i)] The behavior $\B=\ker P(\sigma)=\{
  w(z)\in\F^n[[z,z^{-1}]]\mid P(\sigma)w(z)=0\}$ is controllable.
\item[(ii)] $P(z)$ is left prime over $\F[z,z^{-1}]$.
\item[(iii)] The behavior $\B$ has an image representation. This
  means there exists an $n\times k$ matrix $G(z)$ defined over
  $\F[z,z^{-1}]$ such that
  $$
  \B=\{ \ w(z)\in\F^n[[z,z^{-1}]]\ \mid\ \exists
  m(z)\in\F^k[[z,z^{-1}]]:\ w(z)=G(\sigma)m(z)\ \}.
  $$
\end{description}
\end{thm}

Combining the theorem with the facts that completeness
corresponds to compactness and irreducibility corresponds to
controllability gives a proof of Theorem~\ref{image}.

We conclude the section by defining some parameters of a linear,
time-invariant and complete behavior. For simplicity we will do
this in an algebraic manner. We will first treat behaviors
$\B\subset \F^n[[z]]$, i.e. behaviors in the sense of
Definition~\ref{Def-C-pr}. In Remark~\ref{adjust} we will explain
how the definitions have to be adjusted for behaviors defined on
the time axis $\Z$.

Assume that $P(z)$ is a $r\times n$ polynomial matrix of rank $r$
defining the behavior $\B=\ker P(\sigma)$. There exists a matrix
$U(z)$, unimodular over $\F[z]$, such that
$\tilde{P}(z)=U(z)P(z)$ is row-reduced with ordered row degrees
$\nu_1\geq \ldots\geq \nu_r$.  The indices
$\nu=(\nu_{1},\dots,\nu_{r})$ are invariants of the row module of
$P(z)$ (and hence also invariants of the behavior $\B$), and are
sometimes referred to as the {\em Kronecker indices} or {\em
  observability indices} of $\B$. The invariant
$\delta:=\sum_{i=1}^r\nu_i$ is called the {\em McMillan degree\/}
of the behavior $\B$. If we think of $\B$ as a convolutional code
in the sense of Definition~\ref{Def-C-pr} then we say that $\B$
has transmission rate $\frac{n-r}{n}$. Finally, the free distance
of the code is defined as in\eqr{free-dist}.

\begin{rem}
  The Kronecker indices $\nu$ are in general different from the
  minimal row indices (in the sense of Forney~\cite{fo75}) of the
  $\F(z)$-vector space generated by the rows of $P(z)$. They
  coincide with the minimal row indices if and only if $P(z)$ is
  left prime.
\end{rem}

\begin{rem}                  \label{adjust}
  If $\B\subset \F^n[[z,z^{-1}]]$ is a linear, time-invariant and
  complete behavior, then we can define parameters like the
  Kronecker indices and the McMillan degree in the following way:
  Assume $P(z)$ has the property that $\B=\ker P(\sigma)$. There
  exists a matrix $U(z)$, unimodular over $\F[z,z^{-1}]$, such
  that $\tilde{P}(z)=U(z)P(z)$ is row-reduced and $P(0)$ has full
  row rank $r$. One shows again that the row degrees of
  $\tilde{P}(z)$ are invariants of the behavior. The McMillan
  degree, the transmission rate and the free distance are then
  defined in the same way as for behaviors $\B\subset \F^n[[z]]$.
\end{rem}

\Section{The module point of view}
\label{Module}

Fornasini and Valcher~\cite{fo98a,va94} and the present author in
joint work with Schumacher, Weiner and
York~\cite{ro96a1,ro99a,we98t} proposed a module-theoretic
approach to convolutional codes. The module point of view
simplifies the algebraic treatment of convolutional codes to a
large degree, and this simplification is probably almost necessary
if one wants to study convolutional codes in a multidimensional
setting~\cite{fo98a,va94,we98t}.

{}From a systems theoretic point of view, the module-theoretic
approach studies linear time-invariant systems whose states start
at zero and return to zero in finite time. Such dynamical systems
have been studied by Hinrichsen and
Pr\"{a}tzel-Wolters~\cite{hi80a,hi83a1}, who recognized these
systems as convenient objects for the study of systems
equivalence.

In our development we will again deal with the time axes $\Z$ and
$\Z_+$ in a parallel manner.

\begin{deA}                         \label{Def-D}
  A submodule $\C$ of $\F^n[z,z^{-1}]$ is called a convolutional
  code.
\end{deA}

We like the module-theoretic language. If one prefers to define
everything in terms of trajectories then one could equivalently
define $\C$ as $\F$-linear, time-invariant subset of
$\F^n[[z,z^{-1}]]$ whose elements have finite support.

The analogous definition for codes supported on the positive time
axis $\Z_+$ is:

\renewcommand{\thedefA}{\Alph{defA}$^\prime$} \addtocounter{defA}{-1}
\begin{deA}                          \label{Def-D-pr}
  A submodule $\C$ of $\F^n[z]$ is called a convolutional code.
\end{deA}
\renewcommand{\thedefA}{\Alph{defA}}

Since both the rings $\F[z,z^{-1}]$ and $\F[z]$ are principal
ideal domains (PID), a convolutional code $\C$ has always a
well-defined rank $k$, and there is a full-rank matrix $G(z)$ of
rank $k$ such that $\C=\mathrm{colsp}_{\F[z,z^{-1}]}G(z)$
(respectively $\C=\mathrm{colsp}_{\F[z]}G(z)$ if $\C$ is defined
as in Definition~\ref{Def-D-pr}). We will call $G(z)$ an encoder
of $\C$, and the map
\begin{equation}              \label{map-D}  
\varphi:\ \F^k[z,z^{-1}]\longrightarrow\F^n[z,z^{-1}], \ 
m(z)\longmapsto c(z)=G(z)m(z)
\end{equation}
an encoding map.
\begin{rem}
  In contrast to the situation of Section~\ref{Sym-Dyn}, it is
  possible to define a convolutional code in the sense of
  Definition~\ref{Def-D}
(respectively Definition~\ref{Def-D-pr}) using a rational
encoder. For this, assume that $G(z)$ is an $n\times k$ matrix with 
entries in $\F(z)$. Then
$$
\C=\left\{\ c(z)\in\F^n[z,z^{-1}]\ \mid\ \exists
  m(z)\in\F^k[z,z^{-1}]:\  c(z)=G(z)m(z)\ \right\}
$$
defines a submodule of $\F^n[z,z^{-1}]$. Note that the
map\eqr{map-D} involving a rational encoding matrix $G(z)$ has to
be `input-restricted' in this case.
\end{rem}

In analogy to Lemma~\ref{lem-eq-B} we have:

\begin{lem}                      \label{lem-eq-D}
  Two $n\times k$ matrices $G(z)$ and $\tilde{G}(z)$ defined over
  the Laurent polynomial ring $\F[z,z^{-1}]$ (respectively over
  the polynomial ring $\F[z]$) generate the same code $\C\subset
  \F^n[z,z^{-1}]$ (respectively $\C\subset \F^n[z]$) if and only
  if there is a $k\times k$ matrix $U(z)$, unimodular over
  $\F[z,z^{-1}]$ (respectively over $\F[z]$), such that
  $\tilde{G}(z)=G(z)U(z)$.
\end{lem}

As we already mentioned earlier convolutional codes in the
sense of Definitions~\ref{Def-D} and~\ref{Def-D-pr} are linear
and time-invariant. The following theorem answers any question
about controllability (i.e. irreducibility) and completeness.
\begin{thm}
  A nonzero convolutional code with either Definition~\ref{Def-D}
  or~\ref{Def-D-pr} is controllable and incomplete.
\end{thm}
\begin{proofs}
  The proof of the completeness part of the Theorem is analogous
  to the proof of Theorem~\ref{closure}. In order to show
  controllability, let $G(z)$ be an encoding matrix for a code
  $\C\subset\F^n[z]$ and consider two code words
  $w(z)=G(z)(a_0+a_1+\cdots+a_sz^s)$ and
  $w'(z)=G(z)(b_0+b_1+\cdots+b_sz^s)$. The codeword $w''(z)$
  required by Definition~\ref{control} can be constructed in
  the form
  $$
  G(z)(a_0+a_1+\cdots+a_jz^j+b_{j+\ell}z^{j+\ell}+\cdots+\cdots+b_sz^s).
  $$
\end{proofs}

Submodules of $\F^n[z,z^{-1}]$ (respectively of $\F^n[z]$) form
the Pontryagin dual of linear, time-invariant and complete
behaviors in $\F^n[[z,z^{-1}]]$ (respectively $\F^n[[z]]$).
In the following we follow~\cite{ro96a1} and explain this in a
very explicit way when the time axis is $\Z$. Of course
everything can be done {\em mutatis mutandis} when the time axis is
$\Z_+$.

Consider the bilinear form:
\begin{equation}                       \label{bilin}
\begin{array}{rcl}
(\,,\,) :\hspace{3mm}\F^n[[z,z^{-1}]]\times \F^n[z,z^{-1}]
&\longrightarrow &\F \\
(w,v) &\mapsto &\sum\limits_{i=-\infty}^\infty
\langle w_i,v_i \rangle,
\end{array}
\end{equation}
where $\langle\,,\,\rangle$ represents the standard dot product
on $\F^n$.  One shows that $(\,,\,)$ is well defined and
nondegenerate, in particular because there are only finitely many
nonzero terms in the sum.  For any subset $\C$ of
$\F^n[z,z^{-1}]$ one defines the annihilator
\begin{equation}
\C^\perp = \{w\in \F^n[[z,z^{-1}]]\mid \left( w,v\right)
=0,\forall v\in {\cal C}\}  \label{bdual}
\end{equation}
and the annihilator of a subset $\B$ of $\F^n[[z,z^{-1}]]$ is
\begin{equation}
{\cal B}^{\perp }=\{v\in \F^n[z,z^{-1}]\mid \left( w,v\right)
=0,\forall w\in {\cal B}\}.  \label{bdual2}
\end{equation}

The relation between these two annihilator operations is given
by:

\begin{thm}                       \label{dual}
  If $\C \subseteq \F^n[z,z^{-1}]$ is a convolutional code with
  generator matrix $G(z)$, then $\C^{\perp}$ is a linear,
  left-shift-invariant and complete behavior with kernel
  representation $P(z)=G^t(z)$.  Conversely, if $\B \subseteq
  \F^n[[z,z^{-1}]]$ is a linear, left-shift-invariant and
  complete behavior with kernel representation $P(z)$, then
  $\B^{\perp}$ is a convolutional code with generator matrix
  $G(z)=P^t(z)$.
\end{thm}
\begin{rem}
  An elementary proof of Theorem~\ref{dual} in the case of the
  positive time axis $\Z_+$ is given in~\cite{ro96a1}.
\end{rem}

\begin{rem}
  Theorem~\ref{dual} is a special instance of a broad duality
  theory between solution spaces of difference equations on the
  one hand and modules on the other, for which probably the most
  comprehensive reference is Oberst~\cite{ob90}. In this article
  Oberst~\cite[p. 22]{ob90} works with a bilinear form which is
  different from\eqr{bilin}. This bilinear form induces however
  the same duality as shown in~\cite{gl99u}. Extensions of
  duality results to group codes were derived by Forney and Trott
  in~\cite{fo00u}.
\end{rem}

For finite support convolutional codes in the sense of
Definition~\ref{Def-D} or Definition~\ref{Def-D-pr} the crucial
issue is {\em observability}. In the literature there have been
several definitions of
observability~\cite{fo95p2,fo95p3,fo98a,fo99,lo96a,ro96a1} and it
is not entirely clear how these definitions relate to each other.

In the sequel we will follow~\cite{fo95p2,ro96a1}.
\begin{de} (cf.\,\cite[Prop.\,2.1]{fo95p2})
  A code $\C$ is {\em observable\/} if there exists an integer
  $N$ such that, whenever the supports of $v$ and $v'$ are
  separated by a distance of at least $N$ and $v+v'\in\C$, then
  also $v\in\C$ and $v'\in\C$.
\end{de} 

With this we have the `Pontryagin dual statement' of
Theorem~\ref{wi2}:
\begin{thm}                               \label{dual-wi2}
  (cf.\,\cite[Prop.\,2.10]{ro96a1}) Let $G(z)$ be a $n\times k$
  matrix of rank $k$ defined over $\F[z,z^{-1}]$. The following
  conditions are equivalent:
\begin{description}
\item[(i)] The convolutional code
  $\C=\mathrm{colsp}_{\F[z,z^{-1}]}G(z)$ is observable.
\item[(ii)] $G(z)$ is right prime over $\F[z,z^{-1}]$.
\item[(iii)] The code $\C$ has a kernel representation. This
  means there exists an $r\times n$ `parity-check matrix' $H(z)$
  defined over $\F[z,z^{-1}]$ such that
  $$
  \C=\{ \ v(z)\in\F^n[z,z^{-1}]\ \mid\ \ H(z)v(z)=0\ \}.
  $$
\end{description}
\end{thm}

\begin{rem}
  The concept of observability is clearly connected to the coding
  concept of non-catastrophicity. Indeed an encoder is
  non-catastrophic if and only if the code generated by this
  encoder is observable. In the context of Definition~\ref{Def-A}
  (respectively Definition~\ref{Def-A-pr}) every code has a
  catastrophic as well as a non-catastrophic encoder.  In the
  module setting of Definition~\ref{Def-D} every encoder of an
  observable code is non-catastrophic and every encoder of an
  non-observable code is catastrophic. If one defines a
  convolutional code by Definition~\ref{Def-D} then one could
  talk of a `non-catastrophic convolutional code'. The term
  observable seems however much more appropriate.
\end{rem}

As at the end of Section~\ref{Wil-Beh}, we now define the code
parameters. We do it only for codes given by
Definition~\ref{Def-D-pr} and leave it to the reader to adapt the
definitions to codes given by Definition~\ref{Def-D}.

Assume that $G(z)$ is an $n\times k$ polynomial matrix of rank
$k$ defining the code $\C=\mathrm{colsp}_{\F[z]}G(z)$.  There
exists a unimodular matrix $U(z)$ such that
$\tilde{G}(z)=G(z)U(z)$ is column-reduced with ordered column
degrees $\kappa_1\geq \ldots\geq \kappa_k$.  The indices
$\kappa=(\kappa_{1},\dots,\kappa_{k})$ are invariants of the code
$\C$, which we call the {\em Kronecker indices} or {\em
  controllability indices} of $\C$.  The invariant
$\delta:=\sum_{i=1}^r\kappa_i$ is called the {\em degree\/} of
the code $\C$. The free distance of the code is defined as
in\eqr{free-dist}. Finally we say that $\C$ has transmission rate
$\frac{k}{n}$.

\Section{First-order representations}\label{First-Order}

In this section we provide an overview of the different
first-order representations (realizations) associated with the
convolutional codes and encoding maps which we have defined.

We start with the encoding map\eqr{map-A}. As is customary in
most of the coding literature, we view the map\eqr{map-A} as an
input-output operator from the message space to the code space.
The existence of associated state spaces and realizations can be
shown on an abstract level. Kalman~\cite{ka65,ka69} first showed
how the encoding map\eqr{map-A} can be `factored' resulting in a
realization of the encoding matrix $\varphi$.
Fuhrmamnn~\cite{fu76} refined the realization procedure in an
elegant way. (Compare also~\cite{fu96,ha78}.)

In the sequel we will simply assume that a realization algorithm
exists. We summarize the main results in the following two
theorems:
\begin{thm}                        \label{first-A}
  Let $T(z)$ be a $p\times m$ proper transfer function of
  McMillan degree $\delta$. Then there exist matrices $(A,B,C,D)$
  of size $\delta\times\delta$, $\delta\times m$, $p\times\delta$
  and $p\times m$ respectively such that
  \begin{equation}                 \label{reali1}
    T(z)=C(zI-A)^{-1}B+D.
  \end{equation}
  The minimality conditions are that $(A,B)$ forms a
  controllable pair and $(A,C)$ forms an observable pair.
  Finally\eqr{reali1} is unique in the sense that if
  $T(z)=\tilde{C}(zI-\tilde{A})^{-1}\tilde{B}+\tilde{D}$ with
  $(\tilde{A},\tilde{B})$  controllable and
  $(\tilde{A},\tilde{C})$  observable, then there is a unique
  invertible matrix $S$ such that
   \begin{equation}                 \label{equi1}
   (\tilde{A},\tilde{B},\tilde{C},\tilde{D})=(SAS^{-1},SB,CS^{-1},D).
  \end{equation}
\end{thm}

Consider the encoding map\eqr{map-A} with generator matrix
$G(z)$. Let $m(z)=\sum_{i=s}^tm_iz^i\in{\F}^k((z))$ and
$c(z)=\sum_{i=s}^tc_iz^i\in{\F}^n((z))$ be the sequence of
message and code symbols respectively. Then one has:

\begin{thm}
  Assume that $G(z)$ has the property that $\rank G(0)=k$. Then
  $G(z^{-1})$ is a proper transfer function, and by
  Theorem~\ref{first-A} there exist matrices $(A,B,C,D)$ of
  appropriate sizes such that $G(z^{-1})=C(zI-A)^{-1}B+D$. The
  dynamics of\eqr{map-A} are then equivalently described by:
  \begin{equation}  \label{iso2}
  \begin{array}{rcl}
   x_{t+1} &  = & Ax_t+Bm_t , \\
    c_t &  =  & Cx_t+Dm_t .
  \end{array}
\end{equation}
\end{thm}

The realization\eqr{iso2} is useful if one wants to describe the
dynamics of the encoder $G(z)$. It is however less useful if one
is interested in the construction of codes having certain
properties. The problem is that every code $\C$ has many
equivalent encoders whose realizations appear to be completely
different.
\begin{exmp}
  The encoders
  $$
  G(z)=\left(
    \begin{array}{c}
     \frac{1-z}{z-4}\\  \frac{1+z}{z-4}
    \end{array}
  \right) \mbox{ and } \tilde{G}(z)=\left(
    \begin{array}{c}
     \frac{1-z}{(z-2)(z+3)}\\  \frac{1+z}{(z-2)(z+3)}
    \end{array}
  \right)
  $$
  are equivalent since they define the same code in the sense
  of Definition~\ref{Def-A}. The transfer functions $G(z^{-1})$
  and $\tilde{G}(z^{-1})$ are however very different from a
  systems theory point of view. Indeed, they have different
  McMillan degrees, and over the reals the first is stable
  whereas the second is not. The state space descriptions are
  therefor very different for these encoders.
\end{exmp}
This example should make it clear that for the purpose of
constructing good convolutional codes, representation\eqr{iso2}
is not very useful. \medskip

We are now coming to the realization theory of the behaviors and
codes of Section~\ref{Wil-Beh} and~\ref{Module}. We will continue
with our algebraic approach. The results are stated for the
positive time axis $\Z_+$, but they hold mutatis mutandis for the
time axis $\Z$.

\begin{thm}[Existence]         \label{fothm1}
  Let $P(z)$ be an $r\times n$ matrix of rank $r$ describing a
  behavior $\B$ of the form\eqr{kernelbehav} with McMillan degree
  $\delta$. Let $k=n-r$. Then there exist (constant) matrices
  $G,F$ of size $\delta\times(\delta+k)$ and a matrix $H$ of size
  $n\times(\delta+k)$ such that $\B$ is equivalently described
  by:
  \begin{equation}                    \label{GFH}
   \B=\left\{ w(z)\in\F^n[[z]]\mid \exists 
    \zeta(z)\in\F^{\delta+k}[[z]]:
   \ \ (\sigma G-F)\zeta(z)=0,\ \ w(z)=H\zeta(z)\right\}.
  \end{equation}
  Moreover the following minimality conditions will be satisfied:
\begin{description}
\item[(i)] $G$ has full row rank;
\item[(ii)] $\left[{G \atop H}\right]$ has full column rank;
\item[(iii)] $\left[{zG-F \atop H}\right]$ is right prime.
\end{description}
\end{thm}

For a proof, see~\cite[Thm.\,4.3]{ku94} or~\cite{ku90,ro97a1}.
Equation\eqr{GFH} describes the behavior locally in terms of a
time window of length~1.  The computation of the matrices $G,F,H$
from a kernel description is not difficult. It can even be done
`by inspection', i.e., just by rearranging the
data~\cite{ro97a1}.  The next result describes the extent to
which minimal first-order realizations are unique. A proof is
given in~\cite[Thm.\,4.34]{ku94}.

\begin{thm}[Uniqueness]                 \label{fothm2}
  The matrices $(G,F,H)$ are unique in the following way: If
  $(\tilde{G},\tilde{F},\tilde{H})$ is a second triple of
  matrices describing the behavior $\B$ through\eqr{GFH} and if
  the minimality conditions (i), (ii) and (iii) are satisfied,
  then there exist unique invertible matrices $S$ and $T$ such
  that
\begin{equation}                   \label{sim}
(\tilde{G},\tilde{F},\tilde{G}) \;=\; (SGT^{-1},SFT^{-1},HT^{-1}).
\end{equation}
\end{thm}

The relation to the traditional state-space theory is as follows:
Assume that $P(z)$ can be partitioned into $P(z)=(Y(z)\, U(z))$
with $U(z)$ a square $r\times r$ matrix and $\deg\det
U(z)=\delta$, the McMillan degree of the behavior $\B$. Assume
that $(G,F,H)$ provides a realization for $\B$ through\eqr{GFH}.
Then one shows that the pencil $\left[{zG-F \atop H}\right]$ is
equivalent to the pencil:
\begin{equation}              \label{ABCD1}
  \left[
    \begin{array}{cc}
     zI_\delta-A&B\\ 0 & I_k\\ C & D
    \end{array}
\right].
\end{equation}
The minimality condition (iii) simply translates into the
condition that $(A,C)$ forms an observable pair, showing that the
behavior $\B$ is observable. One also verifies that the matrices
$(A,B,C,D)$ form a realization of the proper transfer function
$U(z)^{-1}Y(z)$ and that this is a minimal realization if and
only if $(A,B)$ forms a controllable pair.  Finally $(A,B)$ is
controllable if and only if the behavior $\B$ is
controllable.\medskip

The Pontryagin dual statements of Theorem~\ref{fothm1}
and~\ref{fothm2} are (see~\cite{ro96a1}):

\begin{thm}[Existence]         \label{fothm3}
  Let $G(z)$ be an $n\times k$ polynomial matrix generating a rate
  $\frac{k}{n}$ convolutional code $\C\subseteq \F^n[z]$ of
  degree $\delta$.  Then there exist $(\delta+n-k)\times \delta$
  matrices $K,L$ and a $(\delta+n-k)\times n$ matrix $M$ (all
  defined over $\F$) such that the code $\C$ is described by
\begin{equation}                    \label{KLM}
\C=\left\{v(z)\in\F^n[z]\mid \exists x(z) \in \F^\delta[z] :\ \
zKx(z)+Lx(z)+Mv(z)=0\right\}. 
\end{equation}
Moreover the following minimality conditions will be satisfied:
\begin{description}
\item[(i)] K has full column rank;
\item[(ii)] $[K \hspace{2mm} M]$ has full row rank;
\item[(iii)] $ [zK+L \ \mid \ M]$ is left prime.
\end{description}
\end{thm}

Equation\eqr{KLM} describes the behavior again locally in terms of a
time window of length~1.

\begin{thm}[Uniqueness]                 \label{fothm4}
  The matrices $(K,L,M)$ are unique in the following way: If
  $(\tilde{K},\tilde{L},\tilde{M})$ is a second triple of
  matrices describing the code $\C$ through\eqr{KLM} and if
  the minimality conditions (i), (ii) and (iii) are satisfied,
  then there exist unique invertible matrices $T$ and $S$ such
  that
\begin{equation}                   \label{sim2}
(\tilde{K},\tilde{L},\tilde{M}) \;=\; (TKS^{-1},TLS^{-1},TM).
\end{equation}
\end{thm}

If $G(z)$ can be partitioned into $G(z)=\left[
    \begin{array}{c}Y(z)\\ U(z)  \end{array}\right]$ with
  $U(z)$ a square $k\times k$ matrix and $\deg\det U(z)=\delta$,
  the degree of the code $\C$, then the pencil $ [zK+L \ \mid \ 
  M]$ is equivalent to the pencil:
\begin{equation}              \label{ABCD2}
  \left[
    \begin{array}{ccc}
     zI_\delta-A& 0_{\delta\times (n-k)} & -B\\ 
             -C & I_{n-k} & -D
    \end{array}
\right].
\end{equation}
The minimality condition (iii) then translates into the condition
that $(A,B)$ forms a controllable pair, showing that the code
$\C$ is controllable. One also verifies that the matrices
$(A,B,C,D)$ form a realization of the proper transfer function
$Y(z)U(z)^{-1}$, that this is a minimal realization if and only
if $(A,C)$ forms an observable pair, and that this is the case if
and only if the code $\C$ is observable.  Finally, the Kronecker
indices of $\C$ coincide with the controllability indices of the
pair $(A,B)$~\cite{ro99a}.

The systems-theoretic meaning of the representation\eqr{ABCD2} is
as follows (see~\cite{ro99a}). Partition the code vector $v(z)$
into:
$$
v(z)=\left[
\begin{array}{c}
   y(z)\\u(z)
\end{array}
\right]\in\F^n[z]
$$
and consider the equation:
\begin{equation}
  \label{kern}
\left[
   \begin{array}{ccc}
   zI_\delta-A&0_{\delta\times (n-k)}&-B\\ -C&I_{n-k}&-D
   \end{array}
 \right]\left[
   \begin{array}{c}
   x(z)\\y(z)\\u(z)
   \end{array}
 \right] =0.
\end{equation}
Let
\begin{eqnarray*}
x(z) &= &x_0z^{\gamma}+x_{1}z^{{\gamma}-1} + \ldots +
x_{\gamma};\ x_t\in\F^\delta, t=0,\ldots,\gamma,\\
u(z) &= &u_0z^{\gamma}+u_{1}z^{{\gamma}-1} + \ldots +
u_{\gamma};\ u_t\in\F^k, t=0,\ldots,\gamma,\\
y(z)& =& y_0z^{\gamma}+y_{1}z^{{\gamma}-1} + \ldots +
y_{\gamma};\ y_t\in\F^{n-k}, t=0,\ldots,\gamma.
\end{eqnarray*}
Then\eqr{kern} is satisfied if and only if
\begin{eqnarray}  \label{iso}
    x_{t+1} & = & Ax_t+Bu_t, \nonumber \\  y_t & = & Cx_t+Du_t, \\
   v_t & = & \left({y_t \atop u_t}\right), 
\,\,  x_0=0,\,\, x_{\gamma+1}=0,\nonumber
\end{eqnarray}
is satisfied. Note that the state-space representation\eqr{iso}
is different from the representation\eqr{iso2}. Equation\eqr{iso}
describes the dynamics of the {\em systematic} and {\em rational}
encoder
$$
G(z)U^{-1}(z)=\left[\begin{array}{c}Y(z)U(z)^{-1}\\
    I_k\end{array} \right].
$$
The encoding map $u(z)\mapsto y(z)=G(z)U^{-1}(z)u(z)$ is
input-restricted, i.e. $u(z)$ must be in the column module of
$U(z)$ in order to make sure that $y(z)$ and $x(z)$ have finite
support.  In terms of systems theory, this simply means that the
state should start at zero and return to zero in finite time.
Linear systems satisfying these requirements have been studied by
Hinrichsen and Pr\"{a}tzel-Wolters~\cite{hi80a,hi83a1}.

\Section{Differences and similarities among the
  definitions}\label{Case1}

After having reviewed these different definitions for
convolutional codes, we would like to make some comparison.

The definitions of Section~\ref{Lin-Alg} and
Section~\ref{Sym-Dyn} viewed convolutional codes as linear,
time-invariant, controllable and observable behaviors, not
necessarily complete.  Definition~\ref{Def-C} and
Definition~\ref{Def-C-pr} were more general in the sense that
non-controllable behaviors were accepted as codes.
Definition~\ref{Def-D} and Definition~\ref{Def-D-pr} were more
general in the sense that non-observable codes were allowed.

In the following subsection we show that all definitions are 
equivalent for all practical purposes if one restricts oneself to 
controllable and observable codes.

\subsection{Controllable and observable codes}

Consider a linear, time-invariant, complete behavior
$\B\subset\F^n[[z,z^{-1}]]$, i.e. a convolutional code in the
sense of Definition~\ref{Def-C}. Let
$$
\C:=\B\cap\F^n((z)).
$$
Then one has
\begin{thm}                      \label{Thm-7-1}
  $\C$ is a convolutional code in the sense of
  Definition~\ref{Def-A-pr}, and its completion $\bar{\C}$ is the
  largest controllable sub-behavior of $\B$.  Moreover, one has a
  bijective correspondence between controllable behaviors
  $\B\subset\F^n[[z,z^{-1}]]$ and convolutional codes
  $\C\subset\F^n((z))$ in the sense of Definition~\ref{Def-A-pr}.
\end{thm}
\begin{proofs}
  Let $\B=\ker P(\sigma)=\{
  w(z)\in\F^n[[z,z^{-1}]]\mid P(\sigma)w(z)=0\}$. If $\B$ is not
  controllable, then $P(z)$ is not left prime and one has a
  factorization $P(z)=V(z)\tilde{P}(z)$, where $\tilde{P}(z)$ is
  left prime and describes the controllable sub-behavior
  $\ker\tilde{P}(\sigma)\subset\B$. Since $\ker V(\sigma)$ is an
  autonomous behavior it follows that 
  $$
  \C=\B\cap\F^n((z))=\ker
  P(\sigma)\cap\F^n((z))=\ker\tilde{P}(\sigma)\cap\F^n((z)).
  $$
  It follows (compare with Theorem~\ref{closure}) that the
  completion $\bar{\C}=\ker\tilde{P}(\sigma)$.
\end{proofs}

Consider now a convolutional code $\C\subset\F^n((z))$ in the
sense of Definition~\ref{Def-A-pr}. Define:
\begin{eqnarray*}
  \check{\C}& := & \C\cap \F^n[z,z^{-1}]\\
\check{\!\!\check{\C}}& := & \C\cap \F^n[z].
\end{eqnarray*}
Conversely if $\C\subset\F^n[z]$ is a convolutional code in the
sense of Definition~\ref{Def-D-pr}, then define:
\begin{eqnarray*}
  \hat{\C}& := & \mathrm{span}_{\F[z,z^{-1}]}\{ v(z)\mid v(z)\in\C\}.\\
\hat{\!\hat{\C}}& := & \mathrm{span}_{\F((z))}\{ v(z)\mid v(z)\in\C\}.
\end{eqnarray*}
By definition it is clear that $\hat{\C}\subset \hat{\!\hat{\C}}$
are convolutional codes in the sense of Definition~\ref{Def-D}
and Definition~\ref{Def-A-pr} respectively.

\begin{thm}                      \label{Thm-7-2}
  Assume that $\C\subset\F^n((z))$ is a convolutional code in the
  sense of Definition~\ref{Def-A-pr}. Then $\ 
  \check{\!\!\check{\C}}\subset \F^n[z]$ is an observable code in
  the sense of Definition~\ref{Def-A-pr}.  Moreover the
  operations $\hat{\hat{}}$ and $\check{\check{}}$ induce a
  bijective correspondence between the observable codes
  $\C\subset\F^n[z]$ and convolutional codes $\C\subset\F^n((z))$
  in the sense of Definition~\ref{Def-A-pr}.
\end{thm}

Theorem~\ref{Thm-7-2} is essentially the Pontryagin dual
statement of Theorem~\ref{Thm-7-1}; we leave it to the reader to
work out the details. Theorem~\ref{Thm-7-1} and~\ref{Thm-7-2}
together show that there is a bijection between controllable and
observable codes in the sense of one definition and another
definition. For controllable and observable codes the code
parameters like the rate $k/n$, the degree $\delta$ and the
Forney (Kronecker) indices are all the same. Moreover the free
distance is in every case the same as well. For all practical
purposes one can therefore say that the frameworks are completely
equivalent, if one is only interested in controllable and
observable codes.

The advantage of Definition~\ref{Def-D} (respectively
Definition~\ref{Def-D-pr}) over the other definitions lies in the
fact that non-observable codes become naturally part of the
theory. It also seems that for construction purposes the relation
between quasi-cyclic codes and convolutional
codes~\cite{ma73,sm99u} is best described in a module-theoretic
framework.

Definition~\ref{Def-C} (respectively Definition~\ref{Def-C-pr})
allows one to introduce non-controllable codes in a natural way.

A Laurent series setting as in Definition~\ref{Def-A-pr} seems to
be most natural if one is interested in the description of the
encoder and/or syndrome former. Extensions of the Laurent series
framework to multidimensional convolutional codes is however much
less natural than the polynomial framework, which is why the
theory of multidimensional convolutional codes has mainly been
developed in a module-theoretic
framework~\cite{fo98a,va94,we98t}.

\subsection{Duality}

In\eqr{bilin} we introduced a bilinear form which induced a
bijection between behaviors $\B\subset\F^n[[z,z^{-1}]]$ and
modules $\C\subset\F^n[z,z^{-1}]$. This duality is a special
instance of Pontryagin duality, and generalizes to group
codes~\cite{fo00u} and multidimensional systems~\cite{ob90}.

In this subsection we show that the bilinear form\eqr{bilin} can
also be used to obtain a duality  between modules and
modules (both in $\F^n[z,z^{-1}]$) or between behaviors and
behaviors (both in $\F^n[[z,z^{-1}]]$).

For this let $\C\subset\F^n[z,z^{-1}]$ be a submodule. Define:
\begin{equation}                        \label{mod-dual}
\C^\vdash :=\C^\perp\cap \F^n[z,z^{-1}].
\end{equation}
One immediately verifies that $\C^\vdash$ is a submodule of
$\F^n[z,z^{-1}]$, which necessarily is observable. One always has
$\C\subset (\C^\vdash)^\vdash$.

One can do something similar for behaviors. For this let
$\B\subset\F^n[[z,z^{-1}]]$ be a behavior. Define:
\begin{equation}                        \label{beh-dual}
\B^\vdash :=\left( \B\cap
  \F^n[z,z^{-1}]\right)^\perp=\overline{\B^\perp}.
\end{equation}
Then it is immediate that $\B^\vdash$ is a controllable
behavior, $(\B^\vdash)^\vdash\subset\B$ and $(\B^\vdash)^\vdash$
describes the controllable sub-behavior of $\B$.

It is also possible to adapt\eqr{bilin} for a duality of
subspaces $\C\subset\F^n((z))$. For such a subspace we define:
\begin{equation}                        \label{Lau-dual}
\C^\vdash :=\left( \C\cap \F^n[z,z^{-1}]\right)^\perp\cap \F^n((z)).
\end{equation}

The duality\eqr{mod-dual} does not in general correspond to the
linear algebra dual of the $R=\F[z,z^{-1}]$ module $\C\subset
R^n$ since there is some `time reversal' involved. The same is
true for the duality\eqr{Lau-dual}, which does not correspond to
the linear algebra dual of the $\F((z))$ vector space $\C$
without time reversal.

If one works however with the `time-reversed' bilinear form:
\begin{equation}                       \label{bilin2}
\begin{array}{rcl}
[\,,\,] :\hspace{3mm}\F^n[[z,z^{-1}]]\times \F^n[z,z^{-1}]
&\longrightarrow &\F \\
(w(z),v(z)) &\mapsto &\sum\limits_{i=-\infty}^\infty
\langle w_i,v_{-i} \rangle
\end{array}
\end{equation}
then the definitions\eqr{mod-dual} and\eqr{Lau-dual} do
correspond to the module dual (and the linear algebra dual
respectively), used widely in the coding literature~\cite{pi88}.
In this case one has: If $G(z)$ is a generator matrix of
$\C^\vdash$ then $H(z):=G^t(z)$ is a parity check matrix of
$(\C^\vdash)^\vdash$.

In the Laurent-series context it is also possible to induce the
duality\eqr{Lau-dual} directly through the time-reversed bilinear
form defined on $\F^n((z))\times \F^n((z))$:
\begin{equation}                       \label{bilin3}
\begin{array}{rcl}
[\,,\,] :\hspace{3mm}\F^n((z))\times \F^n((z))
&\longrightarrow &\F \\
(w(z),v(z)) &\mapsto &\sum\limits_{i=-\infty}^\infty
\langle w_i,v_{-i} \rangle.
\end{array}
\end{equation}
Note that the sum appearing in\eqr{bilin3} is always well
defined. This bilinear form has been widely used in functional
analysis and in systems theory~\cite{fu81a}. 

\subsection{Convolutional codes as subsets of $\F[[z,z^{-1}]]$, a case
  study.}

In this subsection we illustrate the differences of the
definitions in the peculiar case $n=1$.

If one works with Definition~\ref{Def-A} or
Definition~\ref{Def-B} then there exist only the two trivial
codes having the $1\times 1$ generator matrix $(1)$ and $(0)$
as subsets of $\F[[z,z^{-1}]]$.\bigskip

The situation of  Definition~\ref{Def-C} is already more
interesting. For each polynomial $p(z)$ one has the associated
`autonomous behavior':
\begin{equation}                  \label{n=1-Beh}
  \B=\left\{ \ w(z)\ \mid\  p(\sigma)w(z)=0 \ \right\}.
\end{equation}
Autonomous behaviors are the extreme case of uncontrollable
behaviors. If $\deg p(z)=\delta$, then $\B$ is a
finite-dimensional $\F$-vector space of dimension $\delta$. For
coding purposes $\B$ is not useful at all. Indeed, the code
allows only $\delta$ symbols to be chosen freely, say the
symbols $w_0,w_1,\ldots,w_{\delta-1}$. With this the codeword
$w(z)=\sum_{i=-\infty}^\infty w_iz^i\in\B$ is determined, and the
transmission of $w(z)$ requires infinite symbols in the past and
infinite symbols in the future. In other words, the code has
transmission rate $0$. The distance of the code is however very
good, namely $d_\mathrm{free}(\B)=\infty$. If $\B$ is defined on
the positive time axis, i.e. $\B\subset \F[[z]]$ then the
situation is only slightly better. Indeed in this situation, one
sends first $\delta$ message words and then an infinite set of
`check symbols'. As these remarks make clear, a code of the
form\eqr{n=1-Beh} is not very useful.\bigskip

The most interesting situation happens in the setup of
Definition~\ref{Def-D} and Definition~\ref{Def-D-pr}. In this
situation the codes are exactly the ideals $<g(z)>\,
\subset\F[z,z^{-1}]$ (respectively $<g(z)>\, \subset\F[z]$). We
now show that ideals of the form $<g(z)>$ are of interest in the
coding context.
\begin{exmp}
  Let $\F=\F_2=\{0,1\}$. Consider the ideal generated by
  $g(z)=(z+1)$. $<g(z)>\, \subset\F[z,z^{-1}]$ consists in this
  case of the even-weight sequences, namely the set of all
  sequences with a finite and even number of ones. This code is
  controllable but not observable.
\end{exmp}
Ideals of the form $<g(z)>$ are the extreme case of
non-observable behaviors. In principle this makes it impossible
for the receiver to decode a message. However with some
additional `side-information' decoding can still be performed, as
we now explain.\medskip

One of the most often used codes in practice is probably the {\em
  cyclic redundancy check code} (CRC code). These codes are the
main tool to ensure error-free transmissions over the Internet.
They can be defined in the following way: Let $g(z)\in\F[z]$ be a
polynomial. Then the encoding map is simply defined as:
\begin{equation}    \label{enco5}
\varphi:\ \F[z]\longrightarrow\F[z],\ m(z)\longmapsto
c(z)=g(z)m(z). 
\end{equation}
The code is then the ideal $<g(z)>\, =\im(\varphi)$. The distance
of this code is $2$, since there exists an integer $N$ such that
$(z^N-1)\in\, <g(z)>$. As we already mentioned the code is not
observable. Assume now that the sender gives some additional side
information indicating the start and the end of a message. This
can be either done by saying: ``I will send in a moment 1 Mb'',
or it can be done by adding some `stop signal' at the end of the
transmission. Once the receiver knows that the transmission is
over, he applies long division to compute
$$
c(z)=\tilde{m}(z)g(z)+r(z),\ \ \deg r(z)<\delta.
$$
If $r(z)=0$ the receiver accepts the message $\tilde{m}(z)$ as
the transmitted message $m(z)$. Otherwise he will ask for
retransmission. 

The code performs best over a channel (like the Internet) which
has the property that the whole message is transmitted correctly
with probability $p$ and with probability $1-p$ whole blocks of
the message are corrupted during transmission.  One immediately
sees that the probability that a corrupted message $\tilde{m}(z)$
is accepted is~$q^{-\delta}$, where $q=|\F|$ is the field size.

One might argue that the code $<g(z)>\, =\im(\varphi)$ is simply a
cyclic block code, but this is not quite the case. Note that the
protocol does not specify any length of the code word and in each
transmission a different message length can be chosen.  In
particular the code can be even used if the message length is
longer than $N$, where $N$ is the smallest integer such that
$(z^N-1)\in\, <g(z)>$.
\begin{exmp}
  Let $\F=\F_2=\{0,1\}$ and let $g(z)=z^{20}+1$.  Assume
  transmission is done on a channel with very low error
  probability where once in a while a burst error might happen
  destroying a whole sequence of bits. Assume that the sender
  uses a stop signal where he repeats the 4 bits $0011$ for 100
  times.  Under these assumptions the receiver can be reasonably
  sure once a transmission has been complete. The probability of
  failure to detect a burst error is in this case $2^{-20}$ which
  is less than $10^{-6}$. Note that $g(z)$ is a very poor
  generator for a cyclic code of any block length.
\end{exmp}

\begin{rem}
  CRC codes are in practice often implemented in a slightly
  different way than we described it above (see e.g.~\cite{nu90}).
  The sender typically performs long division on $z^\delta m(z)$
  and computes
  $$
  z^\delta m(z)=f(z)g(z)+r(z),\ \ \deg r(z)<\delta.
  $$
  He then transmits the code word $c(z):=z^\delta m(z)-r(z)\in\,
  <g(z)>$. Clearly the schemes are equivalent. The advantage of
  the latter is that the message sequence $m(z)$ is transmitted
  in `plain text', allowing processing of the data immediately.
\end{rem}

\subsection{Some geometric remarks}

One motivation for the author to take a module-theoretic approach
to convolutional coding theory has come from algebraic-geometric
considerations. As is explained in~\cite{lo90,ra94,ra95}, a
submodule of rank $k$ and degree $\delta$ in $\F^n[z]$ describes
a quotient sheaf of rank $k$ and degree $\delta$ over the
projective line ${\mathbb P}^1$. The set of all such quotient
sheaves having rank $k$ and degree at most $\delta$ has the
structure of a smooth projective variety denoted by
$X^{\delta}_{k,n}$. This variety has been of central interest in
the recent algebraic geometry literature. In the context of
coding theory, it has actually been used to predict the existence
of maximum-distance-separable (MDS) convolutional
codes~\cite{ro99a1}.

The set of convolutional codes in the sense of
Definition~\ref{Def-A} or~\ref{Def-A-pr} or~\ref{Def-B} having
rate $\frac{k}{n}$ and degree at most $\delta$ form all proper
Zariski open subsets of $X^{\delta}_{k,n}$. The points in the
closure of these Zariski open sets are exactly the non-observable
codes if the rate is $\frac{1}{n}$. These geometric
considerations suggest that non-observable convolutional codes
should be incorporated into a complete theory of convolutional
codes. The following example will help to clarify these issues:
\begin{exmp}
  Let $\delta=2$, $k=1$ and $n=2$, i.e., consider $X^2_{1,2}$.
  Any code of degree at most 2 then has an encoder of the form:
  $$
  G(z)=\left(
  \begin{array}{c}
    g_1(z)\\g_2(z)
  \end{array}
\right)= \left(
  \begin{array}{c}
   a_0+a_1z+a_2z^2\\ b_0+b_1z+b_2z^2
  \end{array}
\right)
$$
We can identify the encoder through the point
$(a_0,a_1,a_2,b_0,b_1,b_2)\in\mathbb{P}^5$. The variety
$X^2_{1,2}$ is in this example exactly the projective space
$\mathbb{P}^5$.  For codes in the sense of Definition~\ref{Def-A}
or~\ref{Def-A-pr} or~\ref{Def-B}, $G(z)$ must be taken as a basic
minimal encoder in order to have a unique parameterization. This
requires that $g_1(z)$ and $g_2(z)$ are coprime polynomials. The
set of coprime polynomials $g_1(z),g_2(z)$ viewed as a subset of
$\mathbb{P}^5$ forms a Zariski open subset $U\subset
\mathbb{P}^5$ described by the resultant condition
$$
\det \left(
  \begin{array}{cccc}
   a_0 & 0  & b_0 & 0\\
   a_1 & a_0& b_1 & b_0\\
   a_2 & a_1& b_2 & b_1\\
     0 & a_2& 0   & b_2\\
  \end{array}
\right)\neq 0.
$$
For codes in the sense of Definition~\ref{Def-D}, we require
that $a_0$ and $b_0$ are not simultaneously zero in order to have
a unique parameterization.  Definition~\ref{Def-D} leads to a
larger Zariski open set $V$, i.e. $U\subset V\subset
\mathbb{P}^5$. Only with Definition~\ref{Def-D-pr} does one
obtain the whole variety $X^2_{1,2}=\mathbb{P}^5$.
\end{exmp}

In the general situation $X^{\delta}_{k,n}$ naturally contains
the non-observable codes as well. If $k=1$, then
$X^{\delta}_{1,n}=\mathbb{P}^{n(\delta+1)-1}$, and the codes in
the sense of Definition~\ref{Def-D-pr} having rate $\frac{1}{n}$
and degree at most $\delta$ are exactly parameterized by
$X^{\delta}_{1,n}$.

\Section{Conclusion} \label{conc}

The paper surveys a number of different definitions of
convolutional codes.  All definitions have in common that a
convolutional code is a subset $\C\subset\F^n[[z,z^{-1}]]$ which
is both linear and time-invariant. The definitions differ in
requirements such as controllability, observability, completeness
and restriction to finite support.

If one requires that a code be both controllable and observable,
then the restriction to any finite time window will result in
equivalent definitions. Actually Loeliger and
Mittelholzer~\cite{lo96a} define a convolutional code locally in
terms of one trellis section and they require in their definition
that a code is controllable and observable. Algebraically such a
trellis section is simply described through the generalized first
order description\eqr{GFH} or\eqr{KLM}.

If one wants to have a theory which allows one to work with
rational encoders, then it will be necessary that the code has
finite support on the negative time axis $\Z_{-}$ (or
alternatively on the positive time axis $\Z_{+}$). This is one
reason why a large part of the coding literature works with the
field of formal Laurent series.

If one wants in addition to have a theory which can accommodate
non-observable codes (and such a theory seems to have some value)
then it is best to work in a module-theoretic setting.


\end{document}